\documentclass[a4paper]{article}
\usepackage{graphicx} 

\usepackage[utf8]{inputenc} 
\usepackage[english]{babel}
\usepackage{amsfonts}

\usepackage{geometry} 
\usepackage{setspace}
\usepackage{mathtools}
\usepackage{booktabs} 
\usepackage{array} 
\usepackage{paralist} 
\usepackage{verbatim} 
\usepackage{authblk}
\usepackage{subfig} 
\usepackage{amsmath}

\usepackage{fancyhdr} 
\usepackage{url}

\usepackage[nottoc,notlof,notlot]{tocbibind} 
\usepackage[titles,subfigure]{tocloft} 
\usepackage{mathtools}
\usepackage{siunitx}

\usepackage{amsthm,amssymb}

\newcommand*\Laplace{\mathop{}\!\mathbin\bigtriangleup}
\DeclarePairedDelimiterXPP\BigOSI[2]%
  {\mathcal{O}}{(}{)}{}%
  {\SI{#1}{#2}}

\author{Denis Selutckii \\ denis.selutckii@math.msu.ru}

\begin{document}
\title{Upper bounds for Steklov eigenvalues of a hypersurface of revolution}
\affil{Department of Higher Geometry and Topology, Faculty of Mathematics and Mechanics, Moscow State University, Leninskie Gory, GSP-1, 119991, Moscow, Russia}

\date{}
\maketitle
\Large
\begin{abstract}
\Large
    In this paper we find an upper bound for the first Steklov eigenvalue for a surface of revolution with boundary consisting of two spheres of different radii. Moreover, we prove that in some cases this boundary is sharp.
\end{abstract}

{\Large
\section{Introduction}

Steklov eigenvalues are discussed in many recent publications, see e.g. \cite{6, 2, 1, 5}.
Let $\Omega$ be a compact Riemannian manifold with boundary $\Sigma = \partial \Omega$. Then the Dirichlet-Neumann operator is a pseudo-differential operator $\mathcal{D}: C^{\infty}(\Sigma) \rightarrow C^{\infty}(\Sigma)$, which is defined as follows,
\begin{equation*}
    \mathcal{D} f = \partial_{\nu}\hat{f},
\end{equation*}
where $\nu$ is the unit outward normal on $\Sigma$, and $\hat{f}$ is the unique harmonic extension of function $f$ to the interior of $\Omega$. Eigenvalues of $\mathcal{D}$ are called Steklov eigenvalues of manifold $\Omega$. 

It is well known that the Steklov spectrum is discrete, and the eigenvalues form a sequence
\begin{equation*}
    0 = \sigma_0 \leqslant \sigma_1 \leqslant ... \rightarrow +\infty,
\end{equation*}
where all eigenvalues are repeated according to their multiplicity, which is finite. We can also construct the sequence of Steklov eigenvalues without repetition,
\begin{equation*}
    0 = \sigma_{(0)} < \sigma_{(1)} < ... \rightarrow +\infty.
\end{equation*}
Sometimes, we will write eigenvalues of $M$ with a metric $g$ as $\sigma_k(M, g)$.
Also there exists a basis of eigenfunctions  $f_k \in C^{\infty}(\Sigma)$, where $f_k$ is an eigenfunction for the $k$-th eigenvalue $\sigma_k$.

The harmonic extensions $\hat{f_k}$ to $\Omega$ of the functions $f_k$ are solutions of the Steklov spectral problem,
\begin{equation*}
\begin{cases}
    \Laplace{\hat{f_k}} = 0 \quad \text{in } \Omega,\\
    \partial_{\nu}\hat{f_k} = \sigma_k \hat{f_k} \quad \text{in } \Sigma.
\end{cases}
\end{equation*} 

The Rayleigh-Steklov quotient for a function $u$ from the Sobolev space $H^1(\Omega)$ is defined as
\begin{equation*}
    R(u) = \frac{\int_{\Omega}|\nabla u|^2 dV_{\Omega}}{\int_{\Sigma} u^2 dV_{\Sigma}}.
\end{equation*}

Then Steklov eigenvalues are given by the following formula, see e.g. \cite[Section 7.1]{8},
\begin{equation*}
    \sigma_k = \min\{R(f) | f \in H^1(\Omega), \quad f \perp _{\Sigma} f_0, f_1, ..., f_{k-1}\}.
\end{equation*}

Upper bounds for the Steklov eigenvalues for hypersurface of revolution in the Euclidean space with boundary consisting of two components isometric to two unit spheres, were obtained in \cite{1}. In the present paper we generalize these results to the case of hypersurfaces with boundary isometric to two spheres of different radii.

\textbf{Definition 1} (Generalization of the definition from \cite{1}). \textit{An $n$-dimensional compact hypersurface of revolution $(M,g)$ in the Euclidean space with boundary consisting of two components isometric to spheres of radii $R_1$ and $R_2$ is $M = [0, L] \times \mathbb{S}^{n-1}$ with the metric 
\begin{equation*}
    g(r,p) = dr^2 + h^2(r)g_0(p), 
\end{equation*}
where $g_0$ is the canonical metric on $(n-1)$ - dimensional sphere of radius 1, and $h:[0, L] \rightarrow \mathbb{R}_+^*$  is a smooth function, such that
\begin{eqnarray}
    \label{first_cond}
    1) h(0) = R_1, \quad h(L) = R_2, \hspace{1cm}  \\
    \label{second_cond}
    2) |h'(r)| \leqslant 1, \quad \forall r \in [0, L]. \hspace{1cm} 
\end{eqnarray}}

First of all, we prove the following theorem.

\textbf{Theorem 1.} \textit{Let $M = ([0, L] \times \mathbb{S}^{n-1}, g_1=dr^2 + h_1 ^2 g_0)$ be a hypersurface of revolution, $\dim(M) \geq 3$. Then there exists a metric $g_2 = dr^2 + h_2^2 g_0$ on $M$ satisfying \eqref{first_cond} and \eqref{second_cond} such that
\begin{equation*}
    \sigma_k(M,g_1) < \sigma_k(M,g_2), \quad \forall k \geqslant 1.
\end{equation*}
}

Using this fact we can prove later the main result of this paper.

\textbf{Theorem 2.} \textit{Let $(M = [0,L] \times \mathbb{S}^{n-1}, g)$ be a hypersurface of revolution in the Euclidean space of dimension $n \geq 3$ with boundary consisting of two spheres.
Then the first Steklov eigenvalue $\sigma_1 (M, g)$ has the following upper bound,
\begin{multline*}
    \sigma_1 (M, g) < \text{min} \{ \alpha \sigma_1 ^N (A_{R_1+L_1}) + \beta \sigma_1 ^N (A_{R_2+L - L_1}), \\ \frac{1}{1 + (R_1 / R_2)^{n-1}} \sigma_0 ^D (A_{R_1+L_1}) + \frac{1}{1 + (R_2 / R_1)^{n-1}} \sigma_0 ^D (A_{R_2+L - L_1})\},
\end{multline*}
where 
\begin{align*}
    &L_1 = \frac{-R_1+R_2+L}{2}, \quad
    \alpha = \frac{Q_1}{Q_1 + Q_2}, \quad 
    \beta = \frac{Q_2}{Q_1+Q_2}, \\
    &Q_i = R_i^{n-1}(R_i + \frac{1}{(n-1)2^n}(R_1+R_2+L)^n R_i^{1-n})^2, \quad i = 1, 2.
\end{align*}
}

Also, we find an upper bound that does not depend on the length of meridian $L$. 

\textbf{Theorem 3.}\textit{Let $B_n(R_1, R_2)$ be defined as follows 
\begin{eqnarray*}   
    B_n(R_1, R_2):= \frac{1}{1 + (R_1 / R_2)^{n-1}} \sigma_0 ^D (A_{R_1+\frac{L^*-R_1+R_2}{2}}) + \\ + \frac{1}{1 + (R_2 / R_1)^{n-1}} \sigma_0 ^D (A_{R_2+\frac{L^*+R_1-R_2}{2}}),
\end{eqnarray*}
where $L^*$ is the unique solution of the following equation  
\begin{multline*}
    \frac{1}{1 + (R_1 / R_2)^{n-1}} \sigma_0 ^D (A_{R_1+L_1}) + \frac{1}{1 + (R_2 / R_1)^{n-1}} \sigma_0 ^D (A_{R_2+L - L_1}) -  \\ -\alpha \sigma_1 ^N (A_{R_1+L_1}) - \beta \sigma_1 ^N (A_{R_2+L - L_1}) = 0,
\end{multline*}
where $L$ is an unknown variable. Here $L \geqslant R_1 - R_2$ and functions $\alpha$ and $\beta$ are defined in Theorem 2.
}

\textit{Then 
\begin{equation*}
    \sigma_1(M, g) \leqslant B_n(R_1, R_2).
\end{equation*}}

The Steklov eigenfunctions on $(M, g)$ are of the form $u(r) S(p)$, where $u$ is a smooth function, and $S$ is an eigenfunction of the Laplace operator on the sphere. It is an obvious corollary of the standard separation of variables procedure, see e.g. \cite{4}.

The plan of this paper is the following. In Section 2 we prove Theorem 1. In Sections 3 and 4 we give explicit formulas for mixed Steklov-Dirichlet and Steklov-Neumann eigenvalues for spherical shells, respectively. In Section 5 we prove Theorem 2 about an upper bound for the first Steklov eigenvalue, and Theorem 3 about an upper bound that does not depend on $L$.

The author is grateful to Alexei Penskoi for attaching his attention to this problem and for valuable discussions. This research was supported by the Theoretical Physics and Mathematics Advancement Foundation “BASIS” grant Leader (Math) 21-7-1-45-1.

\section{Degenerated maximal metric}

\textbf{Theorem 1.} \textit{Let $M = ([0, L] \times \mathbb{S}^{n-1}, g_1=dr^2 + h_1 ^2 g_0)$ be a hypersurface of revolution, $\dim(M) \geq 3$. Then there exists a metric $g_2 = dr^2 + h_2^2 g_0$ on $M$ satisfying \eqref{first_cond} and \eqref{second_cond} such that
\begin{equation*}
    \sigma_k(M,g_1) < \sigma_k(M,g_2), \quad \forall k \geqslant 1.
\end{equation*}}

\begin{proof}[\unskip\nopunct] \textbf{Proof.}  
A similar theorem in case $R_1 = R_2 = 1$ is proven in the paper \cite{1}. The proof from \cite{1} can be used for our Theorem 1 after some modifications. Without loss of generality let $R_1 \geqslant R_2$. Then $R_1 - R_2 \leqslant L$, because $|h'(r)| \leqslant 1$. \\ 
Let $g_1 = dr^2 + {h^2_1}(r) g_0(p)$ and let $m := \max\limits_{0 \leq r \leq L}{\{h_1(r)\}}$. 
We see that then $R_1 \leqslant m \leqslant (R_1+R_2+L)/2$.

Now we should define a smooth function $h_2(r)$ in the following way.

\begin{enumerate}
    \item If $r \in [0, m - R_1]$ then $h_2(r) = R_1 + r$.
    \item If $r \in [m - R_1, L - m + R_2]$ then we want $h_2$ to be a smooth function greater then $m$.
    \item If $r \in [L - m + R_2, L]$ then $h_2(r) = R_2 + L - r$.
\end{enumerate}


Then for all $r \in [0, 1]$ we have $h_2(r) \geqslant h_1(r)$. Now let $g_2 = dr^2 + {h^2_2}(r) g_0(p)$. The proof continues further exactly as in \cite[proof of Theorem 2]{1}.
\end{proof}

\textbf{Remark 1.} \textit{The same procedure can be repeated for $g_2$. Then we can get a sequence of metrics of revolution $\{g_i\}$. For all $i > j$ we have $\sigma_k(M, g_i) > \sigma_k(M, g_j)$ for all natural $k$. This sequence converges to $\tilde{g} = dr^2 + \tilde{h}^2 g_0$, where \begin{equation*}
\tilde{h} = 
\begin{cases}
    R_1 + r, \quad r \in [0, (L-R_1+R_2)/2], \\
    R_2 + L - r, \quad r \in [(L - R_1 + R_2)/2, L].    
\end{cases}
\end{equation*}
Therefore, the metric $\tilde{g}$ is continious, however it is not smooth. Also $\tilde{g}$ maximizes all Steklov eigenvalues. We call $\tilde{g}$ the degenerated maximal metric for $M$.}

\section{The mixed Steklov-Dirichlet problem}
Let $A_{R+L}$ be a spherical shell $\mathbb{B}^n_{R+L} \setminus \mathbb{\overline{B}}^n_R$, and let $\sigma_{(k)}^ D(A_{R+L})$ denote the $k$-th eigenvalue of mixed Steklov-Dirichlet problem for $A_{R+L}$.

The following proposition is an obvious corollary of the standard separation of variables procedure.

\textbf{Proposition 1.} \textit{Let us solve the Steklov-Dirichlet problem for $A_{R+L}$,
\begin{eqnarray} \label{DirCond}
    \begin{cases}
        \triangle f = 0 \quad \text{in  } A_{R+L},\\
        f = 0 \quad \text{on  } \partial \mathbb{B}_{R+L}, \\
        \partial_{\nu} f = \sigma f \quad \text{on  } \partial \mathbb{B}_R.
    \end{cases}
\end{eqnarray}
Then for any integer number $k \geqslant 0$ the following equality holds,
\begin{equation*}
    \sigma_{(k)}^D (A_{R+L}) = \frac{k+(k+n-2)(R+L)^{2k+n-2}R^{-2k+2-n}}{(R+L)^{2k+n-2}R^{-2k+3-n}-R}.
\end{equation*}}
\begin{proof}[\unskip\nopunct] \textbf{Proof.} 
Let $f(r, p) = u(r) S(p)$. We will find $u(r)$ in the following form
\begin{equation} \label{Dview}
    u(r) = a r^k + br^{-k+2-n},
\end{equation}
since $r^k$ and $r^{-k+2-n}$ form a fundamental system.
$S(p)$ is an eigenfunction for the Laplace operator on sphere, see e.g. \cite{7}.

Then from boundary conditions (\ref{DirCond}) we can find
\begin{equation} \label{DUcond}
    u(R+L) = 0, \quad u'(R) = -\sigma_{(k)}^D u(R).
\end{equation}

The substitution of (\ref{Dview}) in first equality (\ref{DUcond}) gives us the following
\begin{equation*}
    b = -a (R+L)^{2k+n-2}.
\end{equation*}

As eigenfunctions are defined up to a coefficient of proportionality and $a \neq 0$, we obtain 
\begin{equation*}
   u(r) = r^k  - (R+L)^{2k+n-2}r^{-k+2-n}.
\end{equation*}

From (\ref{Dview}) and from second equality (\ref{DUcond}) we can find the $k$-th Steklov-Dirichlet eigenvalue,

\begin{equation*}
    \sigma_{(k)}^D (A_{R+L}) = \frac{k+(k+n-2)(R+L)^{2k+n-2}R^{-2k+2-n}}{(R+L)^{2k+n-2}R^{-2k+3-n}-R}.
\end{equation*}
\end{proof}

The variational description of the Steklov-Dirichlet eigenvalues is well known, see e.g. \cite{6},
\begin{equation*}
    \sigma_k^D = \min\limits_{E\in H_{k+1,0}(A_{R+L})}\max\limits_{0 \neq f \in E}\frac{\int_{A_{R+L}}|\nabla f|^2dV_{A_{R+L}}}{\int_{\partial A_{R,  L}}|f|^2dV_{\partial A_{R+L}}},
\end{equation*}
where $H_{k+1,0}$ is the set of all $(k+1)$-dimensional subspaces in the Sobolev space $H^1_0(A_{R+L})$.

\section{The mixed Steklov-Neumann problem}

Let $\sigma_{(k)} ^ N(A_{R+L})$ denote the $k$-th eigenvalue for the mixed Steklov-Neumann problem on $A_{R+L}$.

The following fact is also an obvious corollary of the standard separation of variables procedure.

\textbf{Proposition 2.} \textit{Consider the mixed Steklov-Neumann problem in $A_{R+L}$,
\begin{eqnarray} \label{NCond}
    \begin{cases}
        \triangle f = 0 \quad \text{in  } A_{R+L},\\
        \partial_{\nu} f = 0 \quad \text{on  } \partial \mathbb{B}_{R+L}, \\
        \partial_{\nu} f = \sigma f \quad \text{on  } \partial \mathbb{B}_R.
    \end{cases}
\end{eqnarray}
Then for any integer number $k \geqslant 1$ the following equality holds,
\begin{equation*}
    \sigma_{(k)} ^N (A_{R+L}) = \frac{k(k + n - 2) (L^{2k+n-2}R^{-k+1-n}-R^{k-1})}{R^k(k+n-2)+kL^{2k+n-2}R^{-k+2-n}}.
\end{equation*}}
\begin{proof}[\unskip\nopunct] \textbf{Proof.} 
Let $f(r, p) = u(r) S(p)$. We will find $u(r)$ in the following form
\begin{equation} \label{Nview}
    u(r) = a r^k + br^{-k+2-n},
\end{equation}
since $r^k$ and $r^{-k+2-n}$ form a fundamental system.
$S(p)$ is an eigenfunction for the Laplace operator on sphere (see for example \cite{7}).

Then from boundary conditions (\ref{NCond}) we can find
\begin{equation} \label{NUcond}
    u'(R+L) = 0, \quad u'(R) = -\sigma_{(k)}^N u(R).
\end{equation}

The substitution of (\ref{Nview}) in first equality (\ref{NUcond}) gives us the following
\begin{equation*}
    b = \frac{ak}{k+n-2}L^{2k-2+n}.
\end{equation*}

As eigenfunctions are defined due to a coefficient of proportionality and $a \neq 0$, we obtain 
\begin{equation*}
    u(r) = r^k  + \frac{k}{k+n-2} L^{2k-2+n} r^{-k+2-n}.
\end{equation*}

From (\ref{Nview}) and from second equality of (\ref{NUcond}) we can find the $k$-th Steklov-Neumann eigenvalue,

\begin{equation*}
    \sigma_{(k)}^N = - \frac{k-k(R+L)^{n}R^{2-2k-n}}{R+\frac{k}{k+n-2}(R+L)^nR^{3-2k-n}}.
\end{equation*}
\end{proof}

The variational description of the Steklov-Neumann eigenvalues is well known, see e.g. \cite{6},
\begin{equation*}
    \sigma_k^N = \min\limits_{E\in H_{k+1}(A_{R+L})}\max\limits_{0 \neq f \in E}\frac{\int_{A_{R+L}}|\nabla f|^2dV_{A_{R+L}}}{\int_{\partial A_{R,  L}}|f|^2dV_{\partial A_{R+L}}},
\end{equation*}
where $H_{k+1}$ is the set of all $(k+1)$-dimensional subspaces in the Sobolev space $H^1(A_{R+L})$.

\section{The upper bound for the first Steklov eigenvalue}
\textbf{Theorem 2.} \textit{Let $(M = [0,L] \times \mathbb{S}^{n-1}, g)$ be a hypersurface of revolution in the Euclidean space of dimension $n \geq 3$ with the boundary consisting of two spheres of different radii.
Then the first Steklov eigenvalue $\sigma_1 (M, g)$ has the following upper bound,
\begin{multline*}
    \sigma_1 (M, g) < \text{min} \{ \alpha \sigma_1 ^N (A_{R_1+L_1}) + \beta \sigma_1 ^N (A_{R_2+L - L_1}), \\ \frac{1}{1 + (R_1 / R_2)^{n-1}} \sigma_0 ^D (A_{R_1+L_1}) + \frac{1}{1 + (R_2 / R_1)^{n-1}} \sigma_0 ^D (A_{R_2+L - L_1})\},
\end{multline*}
where 
\begin{align*}
    &L_1 = \frac{-R_1+R_2+L}{2}, \quad
    \alpha = \frac{Q_1}{Q_1 + Q_2}, \quad 
    \beta = \frac{Q_2}{Q_1+Q_2}, \\
    &Q_i = R_i^{n-1}(R_i + \frac{1}{(n-1)2^n}(R_1+R_2+L)^n R_i^{1-n})^2, \quad i = 1, 2.
\end{align*}}
\begin{proof}[\unskip\nopunct] \textbf{Proof.} 
Let $\phi_1 = f_1(r)  S(p)$ be the eigenfunction of eigenvalue $0$ for the mixed Steklov-Dirichlet problem in $A_{R_1+L_1}$, and let $\phi_2 = f_2(r)  S(p)$ be the first eigenfunction for the mixed Steklov-Neumann problem in $A_{R_2+ L - L_1}$. 

The proof of Proposition 1 implies that the functions $\phi_1$ and $\phi_2$ are of the following form
\begin{align*}
    &\phi_1 = (1 - (L_1 + R_1)^{n-2}(R_1 + r)^{-n+2}) S(p), \\
    &\phi_2 = (1 - (L - L_1 + R_2)^{n-2}(R_2 + r)^{-n+2}) S(p).
\end{align*}

Now we define a new function
\begin{eqnarray}
\Tilde{\phi} = 
\begin{cases}
    \phi_1(r, p), \quad 0 < r < L_1, \\
    \phi_2(L-r, p), \quad L_1 < r < L.
\end{cases}
\end{eqnarray}

The function $\tilde \phi$ is continious. Moreover, if $r = L_1$ then the function $\tilde \phi$ equals 0. Therefore, we can multiply $\phi_1$ by a constant, and obtain
\begin{equation}
    \int_{\mathbb{S}^{n-1}} \phi_1(0,p) dV_\Sigma - \int_{\mathbb{S}^{n-1}} \phi_2(0,p) dV_\Sigma = 0.
\end{equation}

It follows that $R_1^{n-1} \phi_1(0) = R_2^{n-1} \phi_2(0)$.

Therefore, $\tilde \phi$ can be used as a test function for the Rayleigh quotient for first Steklov eigenvalue in $M$. Moreover, we can change the metric $g = dr^2 + h^2(r) g_0$ by the metric $\tilde{g}$ from Remark 1, which is defined as $\tilde{g} = dr^2 + \tilde{h}^2(r) g_0$. Then the Rayleigh quotient will increase. Then we have

\begin{equation*}
\begin{aligned}
    \sigma_1 (M, g) \leq \frac{\int_0 ^L \int_{\mathbb{S}^{n-1}} ((\partial_r\tilde \phi)^2+\frac{1}{\tilde h^2}|\nabla_{\mathbb{S}^{n-1}} \tilde \phi|^2) h^{n-1}dV_{g_0}dr}{\int_\Sigma \tilde \phi^2(r,p) dV_\Sigma} < \\[15pt] < 
    \frac{\int_0 ^{L_1} \int_{\mathbb{S}^{n-1}} ((\partial_r \phi_1)^2+\frac{1}{(R_1+r)^2}|\nabla_{\mathbb{S}^{n-1}} \phi_1|^2) (R_1+r)^{n-1}dV_{g_0}dr}{\int_\Sigma \tilde \phi^2(r,p) dV_\Sigma}+ \\[15pt] +
    \frac{\int_0 ^{L-L_1} \int_{\mathbb{S}^{n-1}}((\partial_r \phi_2)^2+\frac{1}{(R_2+r)^2}|\nabla_{\mathbb{S}^{n-1}} \phi_2|^2) (R_2+r)^{n-1}dV_{g_0}dr}{\int_\Sigma \tilde \phi^2(r,p) dV_\Sigma} = \\[15pt] =\frac{\int_0 ^{L_1} \int_{\mathbb{S}^{n-1}} ((\partial_r \phi_1)^2+\frac{1}{(R_1+r)^2}|\nabla_{\mathbb{S}^{n-1}} \phi_1|^2) (R_1+r)^{n-1}dV_{g_0}dr}{R_1^{n-1}\int_{\mathbb{S}^{n-1}} \phi_1^2(0,p) dV_{g_0} + R_2^{n-1}\int_{\mathbb{S}^{n-1}}\phi_2^2(0,p) dV_{g_0}} + \\[15pt] + \frac{\int_0 ^{L-L_1} \int_{\mathbb{S}^{n-1}} ((\partial_r \phi_2)^2+\frac{1}{(R_2+r)^2}|\nabla_{\mathbb{S}^{n-1}} \phi_2|^2) (R_2+r)^{n-1}dV_{g_0}dr}{R_1^{n-1}\int_{\mathbb{S}^{n-1}} \phi_1^2(0,p) dV_{g_0} + R_2^{n-1}\int_{\mathbb{S}^{n-1}}\phi_2^2(0,p) dV_{g_0}} =
\end{aligned}
\end{equation*}
\begin{equation*}
\begin{aligned}
     = \frac{\int_0 ^{L_1} \int_{\mathbb{S}^{n-1}} ((\partial_r \phi_1)^2+\frac{1}{(R_1+r)^2}|\nabla_{\mathbb{S}^{n-1}} \phi_1|^2) (R_1+r)^{n-1}dV_{g_0}dr}{(1+(\frac{R_1}{R_2})^{n-1})\int_{\mathbb{S}^{n-1}}R_1^{n-1} \phi_1^2(0,p) dV_{g_0}} + \\[15pt]
+
     \frac{\int_0 ^{L-L_1} \int_{\mathbb{S}^{n-1}} ((\partial_r \phi_2)^2+\frac{1}{(R_2+r)^2}|\nabla_{\mathbb{S}^{n-1}} \phi_2|^2) (R_2+r)^{n-1}dV_{g_0}dr}{(1+(\frac{R_2}{R_1})^{n-1})\int_{\mathbb{S}^{n-1}} R_2^{n-1}\phi_2^2(0,p) dV_{g_0}} \leq \\[15pt] \leq
     \frac{1}{1+(\frac{R_1}{R_2})^{n-1}} \sigma_0^D(A_{R_1+L_1}) + \frac{1}{1+(\frac{R_2}{R_1})^{n-1}} \sigma_0^D(A_{R_2 + L - L_1}).
\end{aligned}
\end{equation*}

On another hand, let $\psi_1 = f_1(r)  S(p)$ be an eigenfunction for the first eigenvalue of the mixed Steklov-Neumann problem in $A_{R_1+L_1}$. And also let $\psi_2 = f_2(r)  S(p)$ be an eigenfunction of the first eigenvalue for the mixed Steklov-Neumann problem in $A_{R_2+ L - L_1}$.  

The proof of Proposition 2 implies that the functions $\psi_1$ and $\psi_2$ are of the form
\begin{align*}
    &\psi_1 = ((R_1 + r) - \frac{1}{n-1}(L_1 + R_1)^{n}(R_1 + r)^{1-n}) S(p), \\
    &\psi_2 = ((R_2 + r) - \frac{1}{n-1}(L - L_1 + R_2)^{n}(R_2 + r)^{1-n}) S(p).
\end{align*}

Now let's define a new function
\begin{eqnarray}
\Tilde{\psi} = 
\begin{cases}
    \psi_1(r, p), \quad 0 < r < L_1, \\
    \psi_2(L-r, p), \quad L_1 < r < L.
\end{cases}
\end{eqnarray}

Note, that $\psi_1(L_1) = \psi_2(L-L_1)$. Hence,  $\tilde \psi$ is a continious function, and
\begin{equation}
    \int_{\mathbb{S}^{n-1}} \psi_1(0,p) dV_\Sigma - \int_{\mathbb{S}^{n-1}} \psi_2(0,p) dV_\Sigma = 0,
\end{equation}
The last equality holds because both of these integrals are equal to 0, since $\psi_1, \psi_2$ are eigenfunctions of the first Steklov-Neumann eiegnvalue. 
 
For functions $\psi_1, \psi_2$ the following equality holds,
\begin{equation*}
    \psi_1(0)  \left(R_2 + \frac{1}{n-1}(R_2+L-L_1)^n R_2^{1-n}\right) = \psi_2(0) \left(R_1 + \frac{1}{n-1}(R_1+L_1)^n R_1^{1-n}\right).
\end{equation*}

In the similar way as before we get
\begin{equation*}
\begin{aligned}
    \sigma_1 (M, g) \leq \frac{\int_0 ^L \int_{\mathbb{S}^{n-1}} ((\partial_r\tilde \psi)^2+\frac{1}{\tilde h^2}|\nabla_{\mathbb{S}^{n-1}} \tilde \psi|^2) h^{n-1}dV_{g_0}dr}{\int_\Sigma \tilde \psi^2(r,p) dV_\Sigma} < \\[15pt] < 
    \frac{\int_0 ^{L_1} \int_{\mathbb{S}^{n-1}} ((\partial_r \psi_1)^2+\frac{1}{(R_1+r)^2}|\nabla_{\mathbb{S}^{n-1}} \psi_1|^2) (R_1+r)^{n-1}dV_{g_0}dr}{\int_\Sigma \tilde \psi^2(r,p) dV_\Sigma}+ \\[15pt] +
    \frac{\int_0 ^{L-L_1} \int_{\mathbb{S}^{n-1}}((\partial_r \phi_2)^2+\frac{1}{(R_2+r)^2}|\nabla_{\mathbb{S}^{n-1}} \psi_2|^2) (R_2+r)^{n-1}dV_{g_0}dr}{\int_\Sigma \tilde \psi^2(r,p) dV_\Sigma} = \\[15pt] = \frac{\int_0 ^{L_1} \int_{\mathbb{S}^{n-1}} ((\partial_r \psi_1)^2+\frac{1}{(R_1+r)^2}|\nabla_{\mathbb{S}^{n-1}} \psi_1|^2) (R_1+r)^{n-1}dV_{g_0}dr}{R_1^{n-1}\int_{\mathbb{S}^{n-1}} \psi_1^2(0,p) dV_{g_0} + R_2^{n-1}\int_{\mathbb{S}^{n-1}}\psi_2^2(0,p) dV_{g_0}} + \\[15pt] +
    \frac{\int_0 ^{L-L_1} \int_{\mathbb{S}^{n-1}} ((\partial_r \psi_2)^2+\frac{1}{(R_2+r)^2}|\nabla_{\mathbb{S}^{n-1}} \psi_2|^2) (R_2+r)^{n-1}dV_{g_0}dr}{R_1^{n-1}\int_{\mathbb{S}^{n-1}} \psi_1^2(0,p) dV_{g_0} + R_2^{n-1}\int_{\mathbb{S}^{n-1}}\psi_2^2(0,p) dV_{g_0}} = \\[15pt] =  \frac{\int_0 ^{L_1} \int_{\mathbb{S}^{n-1}} ((\partial_r \psi_1)^2+\frac{1}{(R_1+r)^2}|\nabla_{\mathbb{S}^{n-1}} \psi_1|^2) (R_1+r)^{n-1}dV_{g_0}dr}{(1+\dfrac{R_2^{n-1}(R_2 + \frac{1}{(n-1)2^n}(R_1+R_2+L)^n R_2^{1-n})^2}{R_1^{n-1}(R_1 + \frac{1}{(n-1)2^n}(R_1+R_2+L)^n R_1^{1-n})^2})\int_{\mathbb{S}^{n-1}} R_1^{n-1}\psi_1^2(0,p) dV_{g_0}} + \\[15pt] +\frac{\int_0 ^{L-L_1} \int_{\mathbb{S}^{n-1}} ((\partial_r \psi_2)^2+\frac{1}{(R_2+r)^2}|\nabla_{\mathbb{S}^{n-1}} \psi_2|^2) (R_2+r)^{n-1}dV_{g_0}dr}{(1+\dfrac{R_1^{n-1}(R_1 + \frac{1}{(n-1)2^n}(R_1+R_2+L)^n R_1^{1-n})^2}{R_2^{n-1}(R_2 + \frac{1}{(n-1)2^n}(R_1+R_2+L)^n R_2^{1-n})^2})\int_{\mathbb{S}^{n-1}} R_2^{n-1}\psi_2^2(0,p) dV_{g_0}} \leq \\[15pt] \leq
     \dfrac{1}{1+\dfrac{R_2^{n-1}(R_2 + \frac{1}{(n-1)2^n}(R_1+R_2+L)^n R_2^{1-n})^2}{R_1^{n-1}(R_1 + \frac{1}{(n-1)2^n}(R_1+R_2+L)^n R_1^{1-n})^2}} \sigma_1^N(A_{R_1+L_1}) +
\end{aligned}
\end{equation*}
\begin{equation*}
\begin{aligned}+
 \dfrac{1}{1+\dfrac{R_1^{n-1}(R_1 + \frac{1}{(n-1)2^n}(R_1+R_2+L)^n R_1^{1-n})^2}{R_2^{n-1}(R_2 + \frac{1}{(n-1)2^n}(R_1+R_2+L)^n R_2^{1-n})^2}} \sigma_1^N(A_{R_2 + L - L_1}).
\end{aligned}
\end{equation*}

Now let
\begin{align*}
    &Q_i = R_i^{n-1}(R_i + \frac{1}{(n-1)2^n}(R_1+R_2+L)^n R_i^{1-n})^2, \quad i=1,2,\\
    &\alpha = \frac{Q_1}{Q_1 + Q_2}, 
    \beta = \frac{Q_2}{Q_1+Q_2}.
\end{align*}

Then we have
\begin{equation*}
    \sigma_1(M, g) \leqslant \alpha \sigma_1^N(A_{R_1+L_1}) + \beta \sigma_1^N(A_{R_2 + L - L_1}).
\end{equation*}

It follows that
\begin{multline*}
    \sigma_1 (M, g) < \text{min} \{\alpha \sigma_1 ^N (A_{R_1+L_1}) + \beta \sigma_1 ^N (A_{R_2+L - L_1}), \\ \frac{1}{1 + (R_1 / R_2)^{n-1}} \sigma_0 ^D (A_{R_1+L_1}) + \frac{1}{1 + (R_2 / R_1)^{n-1}} \sigma_0 ^D (A_{R_2+L - L_1}) \}.
\end{multline*}
\end{proof}
$\quad$ \\

\textbf{Proposition 3.} \textit{In the case $R_1 = R_2$ the bound from Theorem 2 is sharp.}
\begin{proof}[\unskip\nopunct] \textbf{Proof.}
Let denote the bound from Theorem 2 by $\mathcal{A}(L)$. Then sharpness of the bound means, that for any $\varepsilon > 0$ there exists such a metric of revolution $g_1$ on $M$ that $\sigma_1(M, g_1) > \mathcal{A}(L) - \varepsilon$. 

Let $\varepsilon > 0$. We will construct on $M$ the metric of revolution $g_{\varepsilon}  = dr^2 + h_{\varepsilon}^2(r) g_0$, where $h_{\varepsilon}(r) = h_{\varepsilon}(L-r)$, and for any $r \in [0, L/2]$ the following inequality holds,
\begin{equation*}
    (R_1+r)^{n-1} - h_{\varepsilon}^{n-1} < \varepsilon^* := \frac{\varepsilon}{\mathcal{A}(L)}.
\end{equation*}

Let $f_1$ be an eigenfunction for $\sigma_1(M, g_{\varepsilon})$. The metric $g_{\varepsilon}$ is symmetric, i.e. $g_{\epsilon}(r, p) = g_{\epsilon}(L-r, p)$. Hence we can choose the function $f_1$ to be symmetric or antisymmetric, i.e. $|f_1(r)| = |f_1(L-r)|$. Then 
\begin{equation*}
\begin{aligned}
   \frac{\int_0 ^{L/2} \int_{\mathbb{S}^{n-1}} ((\partial_r f_1)^2+\frac{1}{(R_1+r)^2}|\nabla_{\mathbb{S}^{n-1}} f_1|^2) (R_1+r)^{n-1}dV_{g_0}dr}{\int_{\mathbb{S}^{n-1}} f_1^2(0,p) dV_{g_0}} - \\ - \frac{\int_0 ^{L/2} \int_{\mathbb{S}^{n-1}} ((\partial_r f_1)^2+\frac{1}{ h_{\varepsilon}^2}|\nabla_{\mathbb{S}^{n-1}} f_1|^2) h_{\varepsilon}^{n-1}dV_{g_0}dr}{\int_{\mathbb{S}^{n-1}} f_1^2(0,p) dV_{g_0}}  < \\
    < \varepsilon^* \frac{\int_0 ^{L/2} \int_{\mathbb{S}^{n-1}} ((\partial_r f_1)^2+|\nabla_{\mathbb{S}^{n-1}} f_1|)dV_{g_0}dr}{\int_{\mathbb{S}^{n-1}} f_1^2(0,p) dV_{g_0}} < \\ <
    \varepsilon^* \sigma_1(M, g_{\varepsilon})  \leqslant \varepsilon \text{,  as  } h_{\varepsilon} > 1 \text{  and  } \varepsilon = \varepsilon^* \mathcal{A}(L).
\end{aligned}
\end{equation*}
Now there are two possible cases.

1) If $f_1 = u_0 S_0(p)$, where $S_0$ is a trivial harmonic function in sphere (i.e. the constant function). Then $f_1$ can be used as a test function for \\ $\sigma_0^D(A_{R_1+L/2})$. Hence, $\sigma_1(M, g_{\varepsilon}) > \sigma_0^D(A_{R_1+L/2}) - \varepsilon $.

2) If $f_1 = u_1 S_1(p)$, where $S_1$ is a harmonic function in sphere for the first eigenvalue. Then $f_1$ can be used as a test function for $\sigma_0^N(A_{R_1+L/2})$. Hence, $\sigma_1(M, g_{\varepsilon}) > \sigma_0^N(A_{R_1+L/2}) - \varepsilon $.

$\quad$

As a result, we have $\sigma_1(M, g_{\varepsilon}) > \mathcal{A}(L) - \varepsilon $. Hence the bound from Theorem 2 is sharp.
\end{proof}
Moreover, we can find an upper bound for the first eigenvalue that does not depend on $L$, i.e. on the length of meridian of the hypersurface of revolution.

\textbf{Theorem 3.}  \textit{Let $B_n(R_1, R_2)$ be defined as follows 
\begin{eqnarray*}   
    B_n(R_1, R_2):= \frac{1}{1 + (R_1 / R_2)^{n-1}} \sigma_0 ^D (A_{R_1+\frac{L^*-R_1+R_2}{2}}) + \\ + \frac{1}{1 + (R_2 / R_1)^{n-1}} \sigma_0 ^D (A_{R_2+\frac{L^*+R_1-R_2}{2}}),
\end{eqnarray*}
where $L^*$ is the unique solution of the following equation  
\begin{multline*}
    \frac{1}{1 + (R_1 / R_2)^{n-1}} \sigma_0 ^D (A_{R_1+L_1}) + \frac{1}{1 + (R_2 / R_1)^{n-1}} \sigma_0 ^D (A_{R_2+L - L_1}) -  \\ -\alpha \sigma_1 ^N (A_{R_1+L_1}) - \beta \sigma_1 ^N (A_{R_2+L - L_1}) = 0
\end{multline*}
where $L$ is an unknown variable. Here $L \geqslant R_1 - R_2$ and functions $\alpha$ and $\beta$ are defined in Theorem 2.}

\textit{Then
\begin{equation*}
    \sigma_1(M, g) \leqslant B_n(R_1, R_2).
\end{equation*}}
\begin{proof}[\unskip\nopunct] \textbf{Proof.}
Let $L_1 = (R_2-R_1+L)/2$. Then let us define the following function
\begin{equation*}
    f_1(L)=\frac{1}{1 + (R_1 / R_2)^{n-1}} \sigma_0 ^D (A_{R_1+L_1}) + \frac{1}{1 + (R_2 / R_1)^{n-1}} \sigma_0 ^D (A_{R_2+L-L_1}).
\end{equation*}
We know from Proposition 1 the explicit formulas for $\sigma_0 ^D (A_{R_1+L_1})$ and $\sigma_0 ^D (A_{R_2+L-L_1})$. It follows that $\sigma_0 ^D (A_{R_1+L_1})$ and $\sigma_0 ^D (A_{R_2+L-L_1})$ are decreasing functions of $L$.

Consider the function
\begin{eqnarray*}
    f_2(L)=\dfrac{1}{1+\frac{R_2^{n-1}(R_2 + \frac{1}{(n-1)2^n}(R_1+R_2+L)^n R_2^{1-n})^2}{R_1^{n-1}(R_1 + \frac{1}{(n-1)2^n}(R_1+R_2+L)^n R_1^{1-n})^2}} \sigma_1^N(A_{R_1+L_1}) + \\[15pt] + \frac{1}{1+\frac{R_1^{n-1}(R_1 + \frac{1}{(n-1)2^n}(R_1+R_2+L)^n R_1^{1-n})^2}{R_2^{n-1}(R_2 + \frac{1}{(n-1)2^n}(R_1+R_2+L)^n R_2^{1-n})^2}} \sigma_1^N(A_{R_2 + L - L_1}).
\end{eqnarray*}

Then the functions $\sigma_1^N(A_{R_1+L_1})$ and $\sigma_1^N(A_{R_2+L - L_1})$ are increasing functions of $L$. Moreover, the derivatives of $\alpha$ and $\beta$ with respect to $L$ are equal by absolute value, but have opposite signs.

Also, $\sigma_1^N(A_{R_2 + L - L_1}) > \sigma_1^N(A_{R_1+L_1})$, and, hence, $f_2(L)$ is a monotone increasing function of $L$.

Therefore, $f_1(R_1-R_2) > f_2(R_1-R_2)$ and there exists only one $L^*$, such that $L^* > R_1 - R_2$, and $f_1(L^*) = f_2(L^*)$. Then $\sigma_1(M, g) < \text{min}(f_1(L), f_2(L)) \leqslant f_1(L^*)$ for all $L \geqslant R_1 - R_2$.
\end{proof}

\newpage
\bibliographystyle{apa}

}
\end{document}